\documentclass[11pt,a4paper]{article}

\usepackage[utf8]{inputenc}

\usepackage{array,amsbsy,amscd,amsfonts,color,amssymb,amstext,amsmath,latexsym,amsthm,amscd
,multicol,graphicx,makeidx,tikz}
\usepackage[english]{babel}
\usepackage{hyperref}

\numberwithin{equation}{section}

\makeindex

\newtheorem{theorem}{Theorem}
\newtheorem{proposition}[theorem]{Proposition}

\newtheorem{lemma}[theorem]{Lemma}

\numberwithin{theorem}{section}

\theoremstyle{definition}
\theoremstyle{definition}\newtheorem{definition}[theorem]{Definition}
\theoremstyle{definition}\newtheorem{remark}[theorem]{Remark}
\theoremstyle{definition}
\theoremstyle{definition}
\theoremstyle{definition}
\theoremstyle{definition}

\def\proof{\trivlist \item[\hskip \labelsep{\bf Proof.\ }]}

\def\proofof [#1] {\noindent {\bf Proof of #1. } }

\def\al #1.{{\mathcal{#1}}}

\newcommand{\A}{\mathcal{A}}
\newcommand{\B}{\mathcal{B}}

\newcommand{\I}{\mathcal{I}}

\renewcommand{\H}{\mathcal{H}}
\newcommand{\M}{\mathcal{M}}

\newcommand{\R}{\mathbb{R}}

\renewcommand{\S}{{S^1}}

\newcommand{\Exp }{{\mathrm {Exp}}}
\newcommand{\Rep }{{\mathrm {Rep}}}
\newcommand{\rev }{{\mathrm {rev}}}
\newcommand{\Out }{{\mathrm {Out}}}
\newcommand{\col }{{\mathrm {col}}}
\newcommand{\Z}{\mathbb{Z}}

\newcommand{\Aut }{{\mathrm {Aut}}}

\newcommand{\bp}{\begin{proof}}
\newcommand{\ep}{\end{proof}}
\newcommand{\bdp}{\begin{dproof}}
\newcommand{\edp}{\end{dproof}}

\newcommand{\Diff}{\operatorname{Diff}(\S)}

\newcommand{\Mob}{\operatorname{PSL(2,\R)}}

\newcommand{\Ad}{\operatorname{Ad }}
\newcommand{\ad}{\operatorname{ad }}

\newcommand{\e}{\varepsilon}

\renewcommand{\a}{\alpha}

\renewcommand{\l}{\lambda}

\newcommand{\ran}{\rangle}
\newcommand{\lan}{\langle}

\newcommand{\End}{\operatorname{End}}
\newcommand{\Hom}{\operatorname{Hom}}

\title{
{\Huge Examples of subfactors from conformal field theory}}

\author{
{\sc Feng Xu}\footnote{Supported in part by an academic senate grant from UCR.}\\
}

\date{}

\begin{document}
\maketitle

\begin{abstract}
Conformal field theory (CFT) in two dimensions provide a rich source
of subfactors. The fact that there are so many subfactors coming
from CFT have led people to conjecture that  perhaps  all finite
depth subfactors are related to CFT. In this paper we examine
classes of subfactors from known CFT. In particular we identify the
so called $3^{\Z_2\times \Z_2}$ subfactor with an intermediate
subfactor from conformal inclusion, and construct new subfactors
from recent work on holomorphic CFT with central charge $24$.

\end{abstract}

\tableofcontents

\section{Introduction}

Subfactor theory provides an entry point into a world of mathematics
and physics containing large parts of conformal field theory,
quantum algebras and low dimensional topology (cf. \cite{J1}) and
references therein). This paper is about  the interactions between
subfactors, algebraic conformal field theory (CFT) which have proved
to be very fruitful lately (cf. \cite{J1, Jp, W}. A number of
general results about cosets, orbifolds and other constructions have
been obtained in the operator algebraic framework (cf. \cite{KLM,
KLX, Xcos, Xorb}) using subfactor techniques. These results have
been conjectured for some time and have resisted all other attempts
so far. \par More recently there have been many new subfactors
constructed by using Planar Algebras (cf. \cite{Jp}) pioneered by V.
Jones. For an incomplete list of references, see \cite{M2221},
\cite{Liu1, Liu2}. Also see \cite{JMS} for a survey. M. Izumi
(cf. \cite{I98, I15}) has constructed a large class of subfactors
generalizing Haagerup subfactors using Cuntz algebras. T. Gannon and
D. Evans (cf. \cite{EG}) have provided evidence suggesting that
Haagerup subfactor may come from CFT.   In \cite{Jren} V. Jones has
devised a renormalization program based on planar algebras as an
attempt to show that all finite depth subfactors are related to CFT,
i.e.,  the double of a finite depth subfactor should be related CFT.
More generally, the program is the following: given a unitary
Modular Tensor Category (MTC) $\M$, (cf. \cite{Tu}), can we
construct a CFT whose representation category is isomorphic to $\M$?

We shall call such a program ``reconstruction program", analogue to a
similar program in higher dimensions by Doplicher-Roberts (cf.
\cite{DR}). M. Bischoff has shown in \cite{Mar1, Mar2} that this can
be done  for  all subfactors with index less than $4$.  In view of
these recent developments, it is natural to examine subfactors from
known CFT. In fact, it is already known that so called $2221$
subfactor are related to subfactor from conformal inclusions (cf.
\cite{M2221}), and it is an interesting question to see if any of these
recently constructed subfactors are related to CFT. We do find a few
more interesting examples in Section \ref{sec:1} after setting up basics in Section \ref{sec:0}.

Another motivation for our work is that it is clear from \cite{Mar1,
Mar2} that holomorphic CFT play an important role in the
reconstruction program. In Section \ref{sec:2} we construct new
subfactors from holomorphic CFT with central charge $24$ based on
recent work.

A major  progress on reconstruction program would be to identify the
origin of Haagerup subfactor in CFT. Despite the evidence in
\cite{EG}, this remains a challenging question. On the other hand,
in Doplicher-Roberts Theorem  a group is constructed first, and then
 a suitable local  net is chosen for the group to act on. In other words
the net is not constructed directly. It takes lots of efforts to
construct conformal net or chiral algebra from MTC which are not
related to groups (see \cite{Mar1, Mar2, R} for recent results), even in
concrete examples such as those in Section \ref{sec:1}. Since
conformal net seems to contain more than MTC, it is possible that a
general reconstruction program may not work. See the end of Section \ref{28obs}
for a possible source of obstructions to the reconstruction program.

I'd like  thank Prof. V. Jones for discussions about his
renormalization program, and Professors C. Dong, P. Grossman, C. Lam
and Z. Liu for useful discussions. This paper is dedicated to the
memory of R. Haag.

\section{Basics of Operator Algebraic Conformal Field
Theory }\label{sec:0}
\subsection{Sectors} Given an infinite factor $M$, the {\it sectors of $M$}  are given
by
$$\text{Sect}(M) = \text{End}(M)/\text{Inn}(M),$$
namely $\text{Sect}(M)$ is the quotient of the semigroup of the
endomorphisms of $M$ modulo the equivalence relation: $\rho,\rho'\in
\text{End}(M),\, \rho\thicksim\rho'$ iff there is a unitary $u\in M$
such that $\rho'(x)=u\rho(x)u^*$ for all $x\in M$.

$\text{Sect}(M)$ is a $^*$-semiring (there are an addition, a
product and an involution $\rho\rightarrow \bar\rho$) equivalent to
the Connes correspondences (bimodules) on $M$ up to unitary
equivalence. If $\rho$ is an element of $\text{End}(M)$ we shall
denote by $[\rho]$ its class in $\text{Sect}(M)$. We define
$\text{Hom}(\rho,\rho')$ between the objects $\rho,\rho'\in \End(M)$
by
\[
\text{Hom}(\rho,\rho')\equiv\{a\in M: a\rho(x)=\rho'(x)a \ \forall
x\in M\}.
\]
We use $\langle  \lambda , \mu \rangle$ to denote the dimension of
$\text{\rm Hom}(\lambda , \mu )$; it can be $\infty$, but it is
finite if $\l,\mu$ have finite index. See \cite{J1} for the
definition of index for type $II_1$ case which initiated the subject
and  \cite{PP} for  the definition of index in general. Also see
Section 2.3 of \cite{KLX} for expositions. $\langle  \lambda , \mu
\rangle$ depends only on $[\lambda ]$ and $[\mu ]$. Moreover we have
if $\nu$ has finite index, then $\langle \nu \lambda , \mu \rangle =
\langle \lambda , \bar \nu \mu \rangle $, $\langle \lambda\nu , \mu
\rangle = \langle \lambda , \mu \bar \nu \rangle $ which follows
from Frobenius duality. $\mu $ is a subsector of $\lambda $ if there
is an isometry $v\in M$ such that $\mu(x)= v^* \lambda(x)v, \forall
x\in M.$ We will also use the following notation: if $\mu $ is a
subsector of $\lambda $, we will write as $\mu \prec \lambda $  or
$\lambda \succ \mu $.  A sector is said to be irreducible if it has
only one subsector.

\subsection{Local nets}

In this section we recall the basic properties enjoyed by the family
of the von Neumann algebras associated with a conformal Quantum
Field Theory on $S^1$ (cf. \cite{GL}). This is an adaption of DHR
analysis (cf. \cite{DHR}) to chiral CFT which is most suitable for our
purposes.\par By an {\it interval}
 we shall always mean an open connected subset $I$
of $S^1$ such that $I$ and the interior $I^\prime $ of its
complement are non-empty.  We shall denote by  ${\I}$ the set of
intervals in $S^1$.

A {\it M\"obius covariant}  net ${\A}$ of von Neumann algebras on
the intervals of $S^1$ is a map
$$
I\rightarrow {\A}(I)
$$
from ${\I}$ to the von Neumann algebras on a Hilbert space ${\H}$
that verifies the following: \vskip .1in \noindent {\bf A. Isotony}.
If $I_1$, $I_2$ are intervals and $I_1 \subset I_2$, then
$$
{\A}(I_1) \subset {\A}(I_2)\, .
$$

\vskip .1in \noindent {\bf B. M\"obius covariance}.  There is a
nontrivial unitary representation $U$ of ${\bold G}$ (the universal
covering group of $PSL(2, {\bold R})$) on ${\H}$ such that
$$
U(g){\A}(I)U(g)^* = {\A}(gI)\, , \qquad g\in {\bold G}, \quad I\in
{\I} \, .
$$

The group $PSL(2, {\bold R})$ is identified with the M\"obius group
of $S^1$, i.e. the group of conformal transformations on the complex
plane that preserve the orientation and leave the unit circle
globally invariant. Therefore ${\bold G}$ has a natural action on
$S^1$. \vskip .1in \noindent {\bf C. Positivity of the energy}.  The
generator of the rotation subgroup $U(R)(\cdot)$ is positive.

Here $R(\vartheta )$ denotes the (lifting to ${\bold G}$ of the)
rotation by an angle $\vartheta $.

\vskip .1in \noindent {\bf D.  Locality}.  If $I_0$, $I$ are
disjoint intervals then ${\A}(I_0)$ and $A(I)$ commute.

The lattice symbol $\vee $ will denote `the von Neumann algebra
generated by'. \vskip .1in \noindent {\bf E. Existence of the
vacuum}.  There exists a unit vector $\Omega $ (vacuum vector) which
is $U({\bold G})$-invariant and cyclic for $\vee _{I\in
{\I}}{\A}(I)$.

\vskip .1in \noindent {\bf F. Uniqueness of the vacuum (or
irreducibility)}.  The only $U({\bold G})$-invariant vectors are the
scalar multiples of $\Omega $.

\vskip .1in \noindent By a {\it conformal net} (or diffeomorphism
covariant net) ${\A}$ we shall mean a M\"{o}bius covariant net such
that the following holds:
\vskip .1in \noindent {\bf G.  Conformal covariance}. There exists a
projective unitary representation $U$ of ${Diff}(S^1)$ on $\H$
extending the unitary representation of $\bold G$  such that for all
$I\in{\I}$ we have
\begin{align*}
U(g)\A(I) U(g)^*\ &=\ \A(gI),\quad g\in {{Diff}}(S^1), \\
U(g)xU(g)^*\ &=\ x,\quad x\in \A (I),\ g\in {{Diff}}(I'),
\end{align*}
where ${{Diff}}(S^1)$ denotes the group of smooth, positively
oriented diffeomorphism of $S^1$ and ${ {Diff}}(I)$ the subgroup of
diffeomorphisms $g$ such that $g(z)=z$ for all $z\in I'$.
\par
Assume  ${\A}$ is  a  M\"{o}bius covariant net. A  M\"{o}bius
covariant {\it representation} $\pi $ of ${\A}$ is a family of
representations $\pi _I$ of the von Neumann algebras ${\A}(I)$,
$I\in {\I}$, on a Hilbert space ${\H}_\pi $ and a unitary
representation $U_\pi $ of the covering group ${\bold G}$ of $PSL(2,
{\bold R})$, with {\it positive energy}, i.e. the generator of the
rotation unitary subgroup has positive generator, such that the
following properties hold:
\begin{align*}
 I\supset \bar I \Rightarrow \pi _{\bar I} \mid _{{\A}(I)}
= \pi _I \quad &\text{\rm (isotony)} \\
\text{\rm ad}U_\pi (g) \cdot \pi _I = \pi _{gI}\cdot \text{\rm
ad}U(g) &\text{\rm (covariance)}\, .
\end{align*}
A unitary equivalence class of  M\"{o}bius covariant representations
of ${\A}$ is called {\it superselection  sector}. \par The
composition of two superselection sectors are  known as Connes's
fusion (cf. \cite{W}). The composition is manifestly unitary and
associative, and this is one of the most important virtues of the
above formulation. The main question is to study all superselection
sectors of  ${\A}$ and their compositions.
\par
Let ${\A}$ be an irreducible conformal net on a
Hilbert space ${\H}$ and let $G$ be a group. Let
$V:G\rightarrow U({\H})$ be a faithful\footnotemark\footnotetext{
If $V:G\rightarrow U({\H})$ is not faithful, we can take $G':= G/ker
V$ and consider $G'$ instead.} unitary representation of $G$ on
${\H}$.
\begin{definition} We say that $G$ acts properly on ${\A}$
if the following conditions are satisfied:\par (1) For each fixed
interval $I$ and each $s\in G$, $\alpha_s (a):=V(s)aV(s^*) \in
{\A}(I), \forall a\in {\A}(I)$; \par (2) For each  $s\in G$,
$V(s)\Omega = \Omega, \forall s\in G$. We will denote by $\Aut(\A)$
all automorphisms of $\A$ which are implemented by proper actions.
\par
\end{definition}
Define ${\A}^G(I):={\B}(I)P_0$ on ${\H}_0$, where ${\H}_0$ is the space of $G$ invariant vectors and
$P_0$ is the projection onto  ${\H}_0$.  The unitary
representation $U$ of ${\bold G}$ on ${\H}$ restricts to an unitary
representation (still denoted by $U$) of ${\bold G}$ on ${\H}_0$.
Then (cf. \cite{Xorb}): \begin{proposition} The map $I\in {\I}\rightarrow {\A}^G(I)$
on $ {\H}_0$ together with the unitary representation (still denoted
by $U$) of ${\bold G}$ on ${\H}_0$ is an irreducible conformal net.
\end{proposition}
We say that ${\A}^G$   is obtained by  {\it orbifold} construction
from ${\A}$.\par \subsection{Complete rationality } We first recall
some definitions from \cite{KLM} . By an interval of the circle we
mean an open connected proper subset of the circle. If $I$ is such
an interval then $I'$ will denote the interior of the complement of
$I$ in the circle. We will denote by ${\I}$ the set of such
intervals. Let $I_1, I_2\in {\I}$. We say that $I_1, I_2$ are
disjoint if $\bar I_1\cap \bar I_2=\emptyset$, where $\bar I$ is the
closure of $I$ in $S^1$. When $I_1, I_2$ are disjoint, $I_1\cup I_2$
is called a 1-disconnected interval in \cite{Xjw}. Denote by
${\I}_2$ the set of unions of disjoint 2 elements in ${\I}$. Let
${\A}$ be an irreducible conformal net. For $E=I_1\cup I_2\in{\I
}_2$, let $I_3\cup I_4$ be the interior of the complement of
$I_1\cup I_2$ in $S^1$ where $I_3, I_4$ are disjoint intervals. Let
$$
{\A}(E):= \A(I_1)\vee \A(I_2), \hat {\A}(E):= (\A(I_3)\vee
\A(I_4))'.
$$ Note that ${\A}(E) \subset \hat {\A}(E)$.
Recall that a net ${\A}$ is {\it split} if ${\A}(I_1)\vee {\A}(I_2)$
is naturally isomorphic to the tensor product of von Neumann
algebras ${\A}(I_1)\otimes {\A}(I_2)$ for any  intervals $I_1,
I_2\in {\I}$ whose closures are disjoint.  ${\A}$ is {\it strongly
additive} if ${\A}(I_1)\vee {\A}(I_2)= {\A}(I)$ where $I_1\cup I_2$
is obtained by removing an interior point from $I$.
\begin{definition}${\A}$ is said to be completely rational,
or $\mu$-rational, if ${\A}$ is split, strongly additive, and the
index $[\hat {\A}(E): {\A}(E)]$ is finite for some $E\in {\I}_2$ .
The value of the index $[\hat {\A}(E): {\A}(E)]$ (it is independent
of $E$ by Prop. 5 of \cite{KLM}) is denoted by $\mu_{{\A}}$ and is
called the $\mu$-index of ${\A}$.  $\A$ is holomorphic if
$\mu_{\A}=1.$
\end{definition}
\par

The following theorem
is proved in \cite{Xorb}:

\begin{theorem}\label{orb} Let ${\A}$ be an irreducible conformal
net and let $G$ be a finite group acting properly on ${\A}$. Suppose
that ${\A}$ is completely  rational. Then:\par (1): ${\A}^G$ is
completely rational and $\mu_{{\A}^G}= |G|^2 \mu_{{\A}}$;
\par (2): There are only a finite number of irreducible covariant
representations of ${\A}^G$ and they give rise to a unitary modular
category as defined in II.5 of \cite{Tu}.
\end{theorem}

For a modular tensor category $\M$, we define
$\Aut(\M)$ to be   the collection of automorphisms of $\M$. One can do
equivariantization of $\M$ with elements of $\Aut(\M)$ (cf.
\cite{EGNO}). We define $\Out(\M)$ to be $\Aut(\M)/ N$  where $N$ is the normal subgroup consisting these automorphisms fixing the isomorphism classes of each simple objects in $\M$. When $\M$ is the representation category
of $\Rep (\A)$ for a complete rational net $\A,$ it may happen that
$\Out(\M)$ contain elements which do not come from $\Aut (\A).$ See
the end of section \ref{lastss} for an example.
\par
Let $\B$ be a  M\"{o}bius (resp.
conformal) net. $\B$ is called a M\"{o}bius (resp. conformal)
extension of $\A$ if there is  a map
\[
I\in\I\to\A(I)\subset \B(I)
\]
that associates to each interval $I\in \I$ a von Neumann subalgebra
$\A(I)$ of $\B(I)$, which is isotonic
\[
\A(I_1)\subset \A(I_2), I_1\subset I_2,
\]
and   M\"{o}bius (resp. diffeomorphism) covariant with respect to
the representation $U$, namely
\[
U(g) \A(I) U(g)^*= \A(g.I)
\] for all $g\in \Mob$ (resp. $g\in \Diff(S^1)$) and $I\in \I$. $\A$ will be called a  M\"{o}bius (resp. conformal)
subnet of $\B.$
\begin{definition}\label{ext}
Let $\A$ be a  M\"{o}bius covariant net. A  M\"{o}bius covariant net
$\B$ on a Hilbert space $\H$ is an extension of $\A$ if there is a
DHR representation  $\pi$ of $\A$ on $\H$ such that $\pi(\A)\subset
\B$ is a   M\"{o}bius subnet. The extension is irreducible if
$\pi(\A(I))'\cap \B(I) = {\mathbb C} $ for some (and hence all)
interval $I$, and is of finite index if $\pi(\A(I))\subset \B(I)$
has finite index for some (and hence all) interval $I$. The index
will be called the index of the inclusion $\pi(\A)\subset \B$ and
will be denoted by $[\B:\A].$  If $\pi$ as representation of
 $\A$ decomposes as $[\pi]= \sum_\lambda m_\lambda[\lambda]$ where
$m_\lambda$  are non-negative  integers and $\lambda$ are
irreducible DHR representations of $\A$, we say that $[\pi]=
\sum_\lambda m_\lambda[\lambda]$ is the spectrum of the extension.
For simplicity we will write  $\pi(\A)\subset \B$ simply as
$\A\subset \B$.
\end{definition}
\begin{lemma}\label{indexab}
If  $\A$ is completely rational, and a   M\"{o}bius covariant net
$\B$ is an irreducible extension of $\A$. Then $\A\subset\B$ has
finite index , $\B$ is completely rational and
$$\mu_\A= \mu_\B [\B:\A]^2.$$
\end{lemma}
\proof
 $\A\subset\B$ has finite
index follows from Prop. 2.3 of \cite{KL}, and the rest follows from
Prop. 24 of \cite{KLM}.
\subsection{Induction}\label{ind}
Let $\B$ be a M\"obius covariant net and $\A$ a subnet. We assume
that $\A$ is strongly additive and $\A\subset \B$ has finite index.
Fix an interval $I_0\in\I$ and  canonical endomorphism (cf.
\cite{LR}) $\gamma$ associated with $\A(I_0)\subset\B(I_0)$. Then we
can choose for each $I\subset\I$ with $I\supset I_0$ a canonical
endomorphism $\gamma_{I}$ of $\B(I)$ into $\A(I)$ in such a way that
the restriction of $\gamma_{I}$ on $\B(I_0)$ is $\gamma_{I_0}$ and $\rho_{I_1}$ is the
identity on $\A(I_1)$ if $I_1\in\I_0$ is disjoint from $I_0$, where
$\rho_{I}\equiv\gamma_{I}$ restricted to $\A(I)$. Given a DHR
endomorphism $\lambda$ of $\A$ localized in $I_0$, the inductions
$\a_{\lambda},\a_{\lambda}^{-}$ of $\lambda$ are the endomorphisms
of $\B(I_0)$ given by
\[
\a_{\lambda}\equiv
\gamma^{-1}\cdot\Ad\e(\lambda,\rho)\cdot\lambda\cdot\gamma \ ,
\a_{\lambda}^{-}\equiv
\gamma^{-1}\cdot\Ad\tilde{\e}(\lambda,\rho)\cdot\lambda\cdot\gamma
\]
where $\e$ (resp. $\tilde\e$) denotes the right braiding (resp. left
braiding) (cf. Cor. 3.2 of \cite{BE1}). In \cite{X98} a slightly
different endomorphism was introduced and the relation between the
two was given in Section 2.1 of \cite{X3m}.

Note that $\Hom( \a_\lambda,\a_\mu)=:\{ x\in \B(I_0) | x
\a_\lambda(y)= \a_\mu(y)x, \forall y\in \B(I_0)\} $ and $\Hom(
\lambda,\mu)= :\{ x\in \A(I_0) | x  \lambda(y)= \mu(y)x, \forall
y\in \A(I_0)\} .$

The following follows from Lemma 3.4 and Th. 3.3 of \cite{X98} (also
cf. \cite{BE1}) :
\begin{theorem}\label{alpha}
(1): $[\lambda]\rightarrow [\a_\lambda], $ are ring
homomorphisms;\par (2) $\lan \a_\lambda,\a_\mu\ran = \lan \lambda
\rho, \mu\ran.$
\end{theorem}

\par
\subsection{Normal inclusions  and tensor equivalence}\label{mirrorextension}
 Let $\B$ be a completely rational net and $\A\subset \B$
be a subnet which is also completely rational.
\begin{definition}\label{coset}
Define a subnet $\widetilde\A\subset \B$ by $\widetilde\A(I):=
\A(I)'\cap \B(I), \forall I\in \I.$
\end{definition}
We note that since $\A$ is completely rational, it is strongly
additive and so we have $\widetilde\A(I)= (\vee_{J\in \I}\A(J))'\cap
\B(I), \forall I\in \I.$ The following lemma then follows directly
from the definition:
\begin{lemma}\label{cosetnet}
The restriction of $\widetilde\A$ on the Hilbert space $ \overline{
\vee_I\widetilde\A(I)\Omega}$ is an irreducible
  M\"{o}bius  covariant net.
\end{lemma}
The net $\widetilde\A$ as in Lemma \ref{cosetnet} will be called the
{\it coset} of  $\A\subset \B$. See \cite{Xcos} for a class of
cosets from Loop groups. \par The following definition generalizes
the definition in Section 3 of \cite{Xcos}:
\begin{definition}\label{cofinite}
 $\A\subset \B$ is called cofinite if the inclusion
$\widetilde\A(I)\vee \A(I) \subset \B(I)$ has finite index for some
interval $I$.
\end{definition}
The following is Prop. 3.4 of \cite{Xmirror2}:
\begin{proposition}\label{rationalc}
Let $\B$ be completely rational, and let $\A\subset\B$ be a
M\"{o}bius subnet which is also completely rational. Then $\A\subset
\B$ is  cofinite if and only if $\tilde\A$ is completely rational.
\end{proposition}
Let  $\B$ be completely rational, and let $\A\subset\B$ be a
M\"{o}bius subnet which is also completely rational. Assume that
$\A\subset \B$ is cofinite. We will use $\sigma_i,\sigma_j,...$
(resp. $\lambda, \mu...$) to label irreducible DHR representations
of $\B$ (resp. $\A$) localized on a fixed interval $I_0$. Since $
\widetilde\A$ is completely rational by Prop. \ref{rationalc},
$\widetilde\A\otimes \A$ is completely rational, and so every
irreducible DHR representation $\sigma_i$ of $\B$, when restricting
to $\widetilde\A\otimes \A$, decomposes as direct sum of
representations of  $ \widetilde\A\otimes \A$ of the form
$(i,\lambda)\otimes \lambda$ by Lemma 27 of \cite{KLM}. Here
$(i,\lambda)$ is a DHR representation of $\widetilde \A$ which may
not be irreducible and we use the  tensor notation
$(i,\lambda)\otimes \lambda$ to represent a DHR representation of $
\widetilde\A\otimes \A$ which is localized on $I_0$ and defined by
$$
(i,\lambda)\otimes \lambda (x_1\otimes x_2)= (i,\lambda)(x_1)\otimes
 \lambda(x_2), \forall x_1\otimes x_2\in \widetilde\A(I_0) \otimes \A(I_0).
$$
We will also identify $\widetilde \A$ and $\A$ as subnets of $
\widetilde\A\otimes \A$ in the natural way. We note that when no
confusion arise, we will use $1$ to denote the vacuum representation
of a net.

\begin{definition}\label{normal}
A  M\"{o}bius subnet $\A\subset\B$ is normal if $\widetilde
\A(I)'\cap \B(I)= \A$ for some I.
\end{definition}

The following is an application of Lemma 2.24 in \cite{Xmirror2}
(also cf. \cite{OS} and \cite{Mar1}):

\begin{theorem}\label{tensore}
Assume $\A$ is completely rational, $\A\subset\B$ is normal and cofinite, and let $\sum_{\lambda\in \Exp}
[(\lambda, \dot{\lambda})]$ be the spectrum of $\A\otimes \widetilde
\A\subset \B.$ Then $\lambda \in \Exp$  are simple objects of a closed fusion
category of $\Rep(\A)$ and there is an equivalence of braided tensor
category $F$ between the subcategory of $\Rep(\A)$
generated by $\lambda \in \Exp$ and  the subcategory of
$\Rep(\widetilde\A)^{\rev}$ generated by $\dot\lambda, \lambda \in
\Exp,$ where the braiding of $\Rep(\widetilde\A)^{\rev}$ is the
mirror image or reversed braiding of $\Rep(\widetilde\A),$ and $F$ maps $\lambda$ to
$\dot\lambda.$
\end{theorem}

\subsection{Subfactors from conformal nets}\label{28obs}  Given a conformal net ${\A},$
There are three general classes of subfactors:\par (i) If $\pi$ is a
covariant representation of ${\A}$, then by locality we have the
following subfactor $\pi_I({\A}(I))\subset \pi_{I'}({\A}(I'))'. $
These are  known as Jones-Wassermann subfactors;\par

(ii) Let $I_1, I_2\in {\I}$. We say that $I_1, I_2$ are disjoint if
$\bar I_1\cap \bar I_2=\emptyset$, where $\bar I$ is the closure of
$I$ in $S^1$. Suppose that  $I_1, I_2$ are disjoint and let $I_3\cup
I_4$ be the interior of the complement of $I_1\cup I_2$ in $S^1$
where $I_3, I_4$ are disjoint intervals. Let
$$
{\A}(E):= {\A}(I_1)\vee {\A}(I_2), \hat {\A}(E):= ({\A}(I_3)\vee
{\A}(I_4))'.
$$ Note that by locality we have subfactor ${\A}(E) \subset \hat {\A}(E).$  These are subfactors analyzed in \cite{Xjw} and \cite{KLM};

\par (iii) If ${\B}\subset {\A}$ is a subnet, then we
have subfactors ${\B(I)}\subset {\A}(I), \forall I,$ and under
certain conditions we get irreducible finite index subfactors, and
in such cases we can induce a representation of ${\B}$ to a soliton
of ${\A}$: i.e., it is only a representation of net ${\A}$
restricted to a punctured circle which is isomorphic to the real
line. By locality such solitons will also give subfactors (cf.
\cite{KLX} and \cite{BEK1}).
\par There are
close relations between the index of subfactors in (i) and (ii).
This is related to the notion of complete rationality in \cite{KLM}.
As for subfactors coming from (iii), a notable class of such
examples come from conformal inclusions and simple current
extensions. For an example, one can construct all subfactors of
index less than $4$ with principal graphs of type $D,E$ from such
constructions. Subfactors induced from simple current extensions also provide examples with
interesting lattice of intermediate subfactors (cf. \cite{Xadv}). \par

Given a rational conformal net $\A$ and any finite group $G$. Assume
that $G\leq S_n$ where $S_n$ is a symmetric group on $n$ letters.
Note that $S_n$ acts on $\A^{\otimes^n}$ by permuting tensors, and
by Th. \ref{orb}  the fixed point algebra $(\A^{\otimes^n})^G$ is
also rational. A particular interesting case is when $G$ is
generated by an cycle, in this case one can relate the chiral
quantities of $(\A^{\otimes^n})^G$  with that of $\A$, and this
leads to many interesting equations (cf. \cite{Ban}, \cite{Xcyc}) due
to the rationality of $(\A^{\otimes^n})^G$ by Th. \ref{orb}.  In
fact P. Bantay proposes  that any MTC $\M$ will have an associated
class of MTCs coming from the orbifolds as in the case of conformal
nets. Following \cite{Ban}, we shall say such MTC verify Orbifold
Covariance Principle. This suggests the following possible
obstructions to reconstruction program: if one can find a MTC which
does not verify Orbifold Covariance Principle, then such a MTC will
not come from CFT. For an example, if one can find a MTC for which
some of those identities in \cite{Xcyc, Ban} are not verified,
then  such a MTC will not come from CFT. However it is not clear
which identity to check, and in fact some properties of MTC such as
certain equations among chiral identities which may follow from
Orbifold Covariance Principle are in fact proved in \cite{Ng}
without constructing the associated class of oribifold MTCs.

\section{Examples of subfactors from extensions and conformal inclusions}\label{sec:1}
 Let $G= SU(n)$. We denote $LG$ the group of
smooth maps $f: S^1 \mapsto G$ under pointwise multiplication. The
diffeomorphism group of the circle $\text{\rm Diff} S^1 $ is
naturally a subgroup of $\text{\rm Aut}(LG)$ with the action given
by reparametrization. In particular the group of rotations
$\text{\rm Rot}S^1 \simeq U(1)$ acts on $LG$. We will be interested
in the projective unitary representation $\pi : LG \rightarrow U(H)$
that are both irreducible and have positive energy. This means that
$\pi $ should extend to $LG\rtimes \text{\rm Rot}\ S^1$ so that
$H=\oplus _{n\geq 0} H(n)$, where the $H(n)$ are the eigenspace for
the action of $\text{\rm Rot}S^1$, i.e., $r_\theta \xi = \exp(i n
\theta)$ for $\theta \in H(n)$ and $\text{\rm dim}\ H(n) < \infty $
with $H(0) \neq 0$. It follows from \cite{PS} that for fixed level
$k$ which is a positive integer, there are only finite number of
such irreducible representations indexed by the finite set
$$
 P_{++}^{k}
= \bigg \{ \lambda \in P \mid \lambda = \sum _{i=1, \cdots , n-1}
\lambda _i \Lambda _i , \lambda _i \geq 0\, , \sum _{i=1, \cdots ,
n-1} \lambda _i \leq k \bigg \}
$$
where $P$ is the weight lattice of $SU(n)$ and $\Lambda _i$ are the
fundamental weights. We will write
$\lambda=(\lambda_1,...,\lambda_{n-1}), \lambda_0= k-\sum_{1\leq
i\leq n-1} \lambda_i$ and refer to $\lambda_0,...,\lambda_{n-1}$ as
components of $\lambda.$

We will use $\Lambda_0$ or simply $1$  to denote the trivial
representation of $SU(n)$.

 For $\lambda , \mu , \nu \in P_{++}^{k}$,
define $N_{\lambda \mu}^\nu  = \sum _{\delta \in P_{++}^{k}
}S_\lambda ^{(\delta)} S_\mu ^{(\delta)} S_\nu
^{(\delta*)}/S_{\Lambda_0}^{(\delta})$ where $S_\lambda ^{(\delta)}$
is given by the Kac-Peterson formula:
$$
S_\lambda ^{(\delta)} = c \sum _{w\in S_n} \varepsilon _w \exp
(iw(\delta) \cdot \lambda 2 \pi /n)
$$
where $\varepsilon _w = \text{\rm det}(w)$ and $c$ is a
normalization constant fixed by the requirement that
$S_\mu^{(\delta)}$ is an orthonormal system. It is shown in
\cite{Kac2} P. 288 that $N_{\lambda \mu}^\nu $ are non-negative
integers. Moreover, define $ Gr(C_k)$ to be the ring whose basis are
elements of $ P_{++}^{k}$ with structure constants $N_{\lambda
\mu}^\nu $.
  The natural involution $*$ on $ P_{++}^{k}$ is
defined by $\lambda \mapsto \lambda ^* =$ the conjugate of $\lambda
$ as representation of $SU(n)$.\par

We shall also denote $S_{\Lambda _0}^{(\Lambda)}$ by $S_1^{(\Lambda
)}$. Define $d_\lambda = \frac {S_1^{(\lambda )}}{S_1^{(\Lambda
_0)}}$. We shall call $(S_\nu ^{(\delta )})$ the $S$-matrix of
$LSU(n)$ at level $k$. \par
 We shall encounter the $\Bbb Z_n$
group of automorphisms of this set of weights, generated by
$$
\sigma : \lambda = (\lambda_1, \lambda_2, \cdots , \lambda_{n-1})
\rightarrow \sigma(\lambda) = ( k -1- \lambda_1 -\cdots
\lambda_{n-1}, \lambda_1, \cdots , \lambda_{n-2}).
$$
We will use $([(\lambda_1, \lambda_2, \cdots , \lambda_{n-1})])$ to denote the orbit of
Define  $\col(\lambda) = \Sigma_i (\lambda_i - 1) i $.
$\col(\lambda)$ will be referred to as the color of $\lambda$. The
central element $ \exp \frac{2\pi i}{n}$ of $SU(n)$ acts on
representation of $SU(n)$ labeled by $\lambda$ as $\exp( \frac{2\pi
i \col(\lambda)}{n})$. The irreducible positive energy
representations of $ L SU(n)$ at level $k$ give rise to an
irreducible conformal net $\A$ (cf. \cite{KLX}) and its covariant
representations. We will use $\lambda=(\lambda_1,...\lambda_{n-1})$
to denote irreducible representations of $\A$ and also the
corresponding endomorphism of $M=\A(I).$

All the sectors $[\lambda]$ with $\lambda$ irreducible generate the
fusion ring of $\A.$ We will use $([(\lambda_1, \lambda_2, \cdots , \lambda_{n-1})])$ to denote the orbit of the sector
$[(\lambda_1, \lambda_2, \cdots , \lambda_{n-1})]$ under the  $\Bbb Z_n$ action above.
\par

\begin{definition}
$v:=(1,0,...,0), v_0:=(1,0,...,0,1) .$ \end{definition}
\subsection{$SU(2)_8$ and the second fish of Bisch-Haagerup}
We will label the irreps of $SU(2)_8$ by half integers $i, 0\leq
i\leq 4.$ Note that the simple current $4$ has integer conformal
dimension, and $\A_{SU(2)_8}$ has a local $\Z_2$ extension $\B$
given by the simple current $4$. By Th. \ref{alpha} we have $[\alpha_2]= [b_1]+[b_2], d_{b_1}=d_{b_2}=
\frac{1+\sqrt{5}}{2},$ and $[b_i^2]=[b_i]+[1], [b_1b_2]=[\alpha_1].
$
 Denote by $\rho$ the index $2$ subfactor for the inclusion
 $\A\subset \B. $
From the fusion rules above  one can easily determine the principal
graph for $\rho b_1$ to be  the second fish of Bisch-Haagerup (cf.
\cite{BH}), and the even vertices are labeled by integers
$0,1,2,3,4$ which are irreps of $\A_{SU(2)_8}$. The double of such a
fusion category was considered in \cite{KE97}. This answers a
question in the introduction of \cite{GI}.
\subsection{Conformal inclusions $SU(n)_{n+2}\subset
SU(\frac{n(n+1)}{2})_1$}

This class of subfactors were discussed in \cite{X98}. We note that
the centralizer algebras $\Hom (\alpha_v^n,\alpha_v^n), n\geq 0$
containing Hecke algebras (cf. \cite{X98}) and additional
generators. The generators and relations for these algebras are
found as in \cite{Liu2} and \cite{Liux}. Another class of conformal
inclusions $SU(n+2)_{n}\subset SU(\frac{(n+2)(n+1)}{2})_1$ are mirror
extensions of $SU(n)_{n+2}\subset SU(\frac{n(n+1)}{2})_1$ and the
corresponding class of subfactors are closely related (\cite{Liux}).
Most of subfactors from such conformal inclusions are not near
group.  Here  we consider a special case $n=4$ to compare with the
near group subfactors in \cite{I15}.

Let $\A_{SU(4)_6}\subset \A_{SU(10)_1}$ be inclusion of nets
correspond to conformal inclusion $SU(4)_6\subset SU(10)_1.$

The fusion graphs for the generators $\alpha_{\Lambda_1},
\alpha_{\Lambda_2}$ can be  determined as in \cite{X98}, and is already
given by different method in \cite{PZ}. We have the following:
$$
[\alpha_{\Lambda_2}\alpha_{\Lambda_2}]= [1]+[\omega
\alpha_{\Lambda_2}] +[\omega^3 \alpha_{\Lambda_2}] + [\omega^{-1}
\alpha_{\Lambda_2}]+[\omega^{-3} \alpha_{\Lambda_2}]
$$
where $[\omega^{10}]=[1],$ $\omega$ is the vector rep of $SU(10)_1.$

In particular if we choose $A=\omega^5 \alpha_{\Lambda_2},$ then
$A$ and $\omega^2=\eta$ generates a fusion subring with
$$
[\eta^5]=[1], [A^2]= [1]+[\eta A] +[\eta^2
A] + [\eta^{3} A]+[\eta^{4}
A]
$$
Such a fusion category is not near-group, but seems to be closely
related to the near-group fusion category in \cite{I15}.\par

We also note that complex conjugation acts on $SU(4)_6\subset
SU(10)_1,$ and we get inclusions $\A_{SU(4)_6}^{\Z_2}\subset
\A_{SU(10)_1}^{\Z_2} $ (cf. \cite{Liux} for more general case). This
is related to   $\Z_2$ equivariantization  of the fusion category above.

\subsection{Conformal inclusions $SU(n)_{n}\subset
Spin(n^2-1)_1$}\label{lasts}
Denote by $\sigma_1$ the vector representation of $Spin (n^2-1)$ and
$v_0$ the adjoint representation of  $SU(n)_n.$  We note that by the
branching rules of $SU(n)_{n}\subset Spin(n^2-1)_1$ we have
$\alpha_{v_0} \succ \sigma_1.$ It follows that $[\sigma_1 \alpha_v]
= [\alpha_v]$ and $\alpha_v$ contains an intermediate subfactor of
index $2$.

We note that the centralizer algebras $\Hom (\alpha_v^n,\alpha_v^n),
n\geq 0$ containing Hecke algebras (cf. \cite{X98}) and additional
generators. It is an interesting questions to find generators and
relations for these algebras as in \cite{Liu2} and \cite{Liux}.

The case when $n=4$ is analyzed in \cite{X98}. Here we consider the
case $n=5$.

The branching rules for  the conformal inclusion $SU(5)_5\subset Spin(24)_1$
are given by (cf. \cite{KW}):
$$
[1]=([(0,0,0,0)])+ ([(0,1,0,2)]), [\sigma_1]=[(1,1,1,1)]+
([(1,0,0,1)]),  [\sigma_2]= [\sigma_3]=2[(1,1,1,1)]
$$
where $\sigma_1, \sigma_2, \sigma_3$ denote the vector and spinor
representations which form $\Z_2\times \Z_2$ under fusion,
and $([(\lambda_1, \lambda_2, \lambda_3, \lambda_4)])$ denotes the orbit of sector
$[(\lambda_1, \lambda_2, \lambda_3, \lambda_4)]$ under the center $\Z_4$ of $SU(4)$.  Here by slightly abusing notations the left hand side of the
equations above are understood as the restrictions from representations of $Spin(24)_1$ to $SU(5)_5.$
The
fusion of the adjoint representation is given by (cf. \cite{Xadv})
$$
[v_0^2]= [1]+2[v_0]+[(2,0,0,2)]+[(0,1,1,0)]+ [(0,1,0,2)]+[(2,0,1,0)].
$$

By using the above and Th. \ref{alpha} we have
$$
\langle\alpha_{v_0},\alpha_{v_0}\rangle= 3.
$$

It follows that $[\alpha_{v_0}]=[\sigma_1]+ [A] +[\sigma_1 A] $
where $A$  is irreducible. Since $[\sigma_1 \alpha_v]=[\alpha_v],$ $\alpha_{v}$ has an
intermediate subfactor denoted by $\rho_1$ and
$[\rho_1\bar{\rho_1}]=[1]+[A].$

We have the following fusion rules:

$$
[A^2]=[1]+[A\sigma_1]+[A\sigma_1]+[A\sigma_2]+[A\sigma_3]
$$

Let $\A_1$ be the simple current
extensions of $\A_{SU(5)_5}.$ Note that $\A_{SU(5)_5}\subset
\A_1\subset \A_{Spin(24)_1}.$ Note that color $0$ irreducible
representations of $\A_{SU(5)_5}$ induce to DHR representations of
$\A_1.$ We enumerate these $10$ irreducible representations of
$\A_1$ as follows: $1, z_1, z_2, z_3, \ad, b_i, 1\leq i \leq 5$
where $1$ is the vacuum representation, $z_1,z_2,z_3$ are induced
from $(0,0,2,1), (0,1,0,2)$ and $(0,1,1,0)$ respectively , $\ad$ is
induced from $(1,0,0,1),$ and $b_i, 1\leq i \leq 5$ are irreducible
components of the representation induced from $(1,1,1,1).$ Let $g\in
\Z_5$ be the generator of $\Aut (\A_1)$ due to the simple current
extension. We have $[\Ad_g b_i] = [b_{i+1}], 1\leq i \leq 5$ and $\Ad_g$
fix the rest $5$ irreducible representations of $\A_1.$  We consider
the induction from $\A_1$ to $\A_{Spin(24)_1}.$ The branching rules  of the
 inclusion $\A_1\subset \A_{Spin(24)_1}$ are given by :
$$
[\Lambda_0]=[1]+ [z_2], [\sigma_1]=[b_1]+ [\ad],
[\sigma_2]=[b_2]+[b_3], [\sigma_3]= [b_4]+[b_5]
$$

The following can be determined from the fusion rules and
Th. \ref{alpha}:
\begin{align*}
& [\alpha_{\ad}]= [\sigma_1]+[ A]+[\sigma_1 A],
 [\alpha_{b_1}]=
[\sigma_1]+[\sigma_2 A]+[\sigma_3 A], \\
& [\alpha_{b_2}]= [\sigma_2]+[A]+[\sigma_2 A],  [\alpha_{b_3}]=
[\sigma_2]+[\sigma_1 A]+[\sigma_3 A], \\
& [\alpha_{b_4}]= [\sigma_3]+[ A]+[\sigma_3 A], [\alpha_{b_5}]=
[\sigma_3]+[\sigma_1 A]+[\sigma_2 A] \\
& [\alpha_{z_1}]=[\alpha_{z_3}]=[A]+[A\sigma_1]+[A\sigma_2]+[A\sigma_3], \\
& [\alpha_{z_2}]=[1]+ [A]+[A\sigma_1]+[A\sigma_2]+[A\sigma_3]
\end{align*}

We see the intermediate subfactor,  denoted by $\rho_1$ above, is exactly $3^{\Z_2\times \Z_2}$
subfactor constructed by Izumi (cf. \cite{GI}). The  double of this
subfactor is computed in \cite{GI}. By \cite{BEKd} we can now see
the double is $\Rep{\A_1}\otimes \Rep{\B}^{\rev}.$  Consider
inclusions $\B\otimes \B \subset \A_{Spin (48)_1}\subset \B_1$ where
$\B_1$ is $\Z_2$ extension of $\A_{Spin (48)_1}$ which is holomorphic. Inspecting the
spectrum of $\B\otimes \B \subset \B_1$ we see that the inclusion
$\B\subset \B_1$ is normal, and by Th. \ref{tensore} we conclude
that $\Rep{\B}^{\rev}\simeq \Rep{\B}$ as braided tensor categories.
So we have shown the following theorem:

\begin{theorem}
The double of  $3^{\Z_2\times \Z_2}$ subfactor is the category of
representations of $\A_1\otimes \A_{Spin(24)_1}$ and verifies the Orbifold
Covariance Principle.
\end{theorem}


In \cite{GI}, a $\Z_3$ equivarization of $3^{\Z_2\times \Z_2}$ is
found. However there is no such $\Z_3$ in $\Aut(\A_1)$ which lifts
to $\Z_3$ on $\A_{Spin(24)_1}$ permuting $\sigma_i.$  If
reconstruction program works, there must be a conformal net $\B_1$
with $\Z_3\in \Aut(\B_1)$ such that the category of representations of
$\B_1^{\Z_3}$ is braided equivalent to the category of representations
of $\A_1\otimes \A_{Spin(24)_1}.$ It is also an interesting question to see if one
can relate $\A_1$ to the double of $D2D$ subfactor in \cite{GI}.

\subsection{Conformal inclusions $SU(n)_m \times SU(m)_n\subset Spin(nm)_1$}
Such cases were studied in \cite{X981} and \cite{OS}. In particular,
by applying Th. \ref{tensore}, we get a braided tensor equivalence
between the fusion subcategory of $\Rep \A_{SU(n)_m}$ generated by
the color $0$ irreducible representations and  the fusion
subcategory of $\Rep \A_{SU(m)_n}^{\rev}$ generated by the color $0$
irreducible representations.

\section{ Subfactors from holomorphic CFT of central charge $24$}\label{sec:2}
In \cite{Sch} A. Schellekens gave a conjectured list of $71$
holomorphic CFT of central charge $24.$ Many examples on this list
have been realized (cf. \cite{Lam1, Lam2, Eho}, \cite{Xmirror2}). We
note that the examples in \cite{Lam1, Lam2, Eho} are given in the
language of lattice VOAs, their orbifolds, and affine Kac-Moody
algebras.  To translate these results into conformal nets, we need
to know that irreducible representations of lattice VOAs give rise
to irreducible representations of corresponding lattice conformal
nets, and irreducible twisted  representations of lattice VOAs give
rise to irreducible solitons of corresponding lattice conformal
nets. These are proved as Prop. 3.15 and Prop. 4.25 of \cite{DX},
and all these representations are of finite index (cf. Cor. 4.31
\cite{DX}). The fact that irreducible representations of affine
Kac-Moody algebra give rise to irreducible representation of the
corresponding conformal net can be seen from the theorem on P. 488
of \cite{W}.

One of the special property of such holomorphic CFT $\B$ is that if
the weight $1$ subspace is non-zero, then (cf. \cite{DM}) it
generates a Kac-Moody subnet $\A\subset \B$ such that the spectrum
of $\A\subset \B$ is finite.  Hence just like Section \ref{sec:1} we
can consider subfactors associated with such $\A\subset \B.$ To
apply the results of  Section \ref{sec:1} we will also need that the
index $[\B:\A]< \infty.$ This is expected to be true in general but
the general case is only known for type $A$ Kac-Moody algebras and
type $D$ at odd level (cf. \cite{W, TL} ). We select a few examples
from \cite{Sch} which have been constructed recently, in a similar
order as in the previous section.
\subsection{No. 20 in \cite{Sch}}
This case is constructed in \cite{Lam2}. We have $[\B:\A]< \infty$
by \cite{TL}, and so we have a new finite index subfactor. If we
take the commutant of the subnet generated by $SU(2)_1^2$, then we
see that we get a local extension of $\A_{Spin({12})_5}\subset \B_1$
whose spectrum is given by
$$
[1] + [010002] + [010020] +[100111] + [002000] + [200100].
$$

\begin{theorem}
There is a local  extension  $\A_{Spin({12})_5}\subset \B_1$ whose
spectrum is given by
$$
[1] + [010002] + [010020] +[100111] + [002000] + [200100]
$$
where we have used the same notation of \cite{Sch} for representations of $Spin({12})_5.$
\end{theorem}

We can now consider the induced subfactor $\alpha_v$ where $v$
denotes the vector representation of $\A_{Spin({12})_5}.$ The
centralizer algebras $\Hom (\alpha_v^n, \alpha_v^n), n\geq 0$ will
contain BMW algebra as in \cite{X97}. We know that $d_v= 7.7396813.$
It is an interesting question to analyze the nature of such
algebras.

\begin{remark}
The above  local extension of $\A_{Spin({12})_5}\subset \B_1$ should
be the mirror extension (cf. \cite{Xmirror1}) associated with the
conformal inclusion $Spin(5)_{12}\subset E_8.$  This example is
similar to No. 27 on Schelleken's list, which is also constructed by
using the mirror extensions of conformal inclusion $SU(3)_9\subset
E_6$ in \cite{Xmirror2}.
\end{remark}

\subsection{No. 11 in \cite{Sch}}
This case is constructed by C. Lam by using the ideas of
\cite{Lam2}.  We have $[\B:\A]< \infty$ by \cite{W}, and so we have
a new finite index subfactor.
\begin{theorem} There is a local
extension $\A_{SU(7)_7}\subset \B$ whose spectrum is given by
$$
([1]) + ([(0,0,1,3,0,1)]) + ([(1,0,3,1,0,0)]) + ([(0,0,2,0,3,0)]) + ([(0,1,0,1,2,2)]) + 3 [(1,1,1,1,1,1)]
$$
\end{theorem}
We can now consider the induced subfactor $\alpha_v$ where $v$
denotes the vector representation of $SU({7})_7$  and $d_v=
\frac{1}{\sin(\frac{\pi}{14})}= 4.493959. $ The centralizer algebras
$\Hom (\alpha_v^n, \alpha_v^n), n\geq 0$  contain Hecke algebras
as in \cite{X98}. It is an interesting question to analyze the
nature of such algebras.

\subsection{No. 9 in \cite{Sch}}\label{lastss}
This case is constructed in \cite{Lam1} and \cite{Eho}. By examining
the spectrum, we can see that we have $\A_1\otimes \A_1\subset \B$
and the spectrum is given by
$$
[1]\otimes [1]+ [z_2]\otimes [z_2]+ [z_1]\otimes [z_3] +
[z_3]\otimes [z_1] + [b_1]\otimes [\ad] +  [\ad]\otimes
[b_{\tau(1)}] + [b_2]\otimes [b_{\tau(2)}] + [b_3]\otimes
[b_{\tau(3)}]+ [b_4]\otimes [b_{\tau(4)}] + [b_5]\otimes
[b_{\tau(5)}],
$$
where $\A_1$ and its irreducible representations are as in section
\ref{lasts}, and $\tau\in S_5$.  So the inclusion $\A_1 \subset \B$
is normal, and we have a braided tensor category equivalence $F_1:
\Rep (\A_1) \rightarrow \Rep (A_1)^{\rev}$ such that $F_1(z_1)=z_3,
F_1(z_3)=z_1, F_1(b_1)=\ad$  .

Now consider the conformal inclusions $SU(5)_5\times SU(5)_5\subset
SU(25)_1\subset \B_2$ where $\B_2$ denotes the holomorphic net
corresponding to No. 67 on Schelleken's list. By examining the
spectrum we find that we have $\A_1\otimes \A_1\subset \B_2$ where
the spectrum is given by
$$
[1]\otimes [1]+ [z_2]\otimes [z_2]+ [z_1]\otimes [z_3] +
[z_3]\otimes [z_1] + [\ad]\otimes [\ad] +  [b_1]\otimes [b_1] +
[b_2]\otimes [b_2] + [b_3]\otimes [b_3]+ [b_4]\otimes [b_4] +
[b_5]\otimes [b_5].
$$
So the inclusion $\A_1 \subset \B_2$ is normal, and we have a
braided tensor category equivalence $F_2: \Rep (\A_1) \rightarrow
\Rep (\A_1)^{\rev}$ such that $F_2(z_1)=z_3, F_2(z_3)=z_1 $ . We
note that the generator $g\in \Aut(\A_1)$ induces braided tensor
category equivalence of $\Rep (\A_1)$  of order  $5.$ Composing $g$
with $F_1F_2^{-1}$, and examining actions on the $6$ element set of
irreducible representations $\ad, b_i, 1\leq i\leq 5$,  we see that
these two elements generate a subgroup of $S_6$ which acts
transitively on $6$ letters. Such a group has at least order $60$.
 We have therefore
proved the following:
\begin{theorem}
 $\Out(\Rep (\A_1))$
has at least order $60$.
\end{theorem}

It remains an interesting question to determine the equivarizations
of  $\Rep (\A_1)$ with respect to the group elements in the above
theorem. Except for the case of $g\in \Aut(\A_1),$ where the
orbifold net is simply $\A_{SU(5)_5},$ the rest of the cases are not
known to be related to CFT, and this includes the $\Z_3$ case
considered in section \ref{lasts}.

We note that as soon as central charge is greater than $24$, there
are a lot more holomorphic CFT since there are many  more unimodular
even positive definite  lattices.  It is an interesting question to
to see if one can get more interesting subfactors related to
holomorphic CFT with central charge greater than $24$.



\author{
{\small Department of Mathematics}\\
{\small University of California at Riverside}\\
{\small E-mail: {\tt xufeng@math.ucr.edu}} }

\date{}
\end{document}